\def\sqr#1#2{{\vcenter{\hrule height.#2pt
      \hbox{\vrule width.#2pt height#1pt \kern#1pt
         \vrule width.#2pt}
       \hrule height.#2pt}}}
\def\eop{\mathchoice\sqr34\sqr34\sqr{2.1}3\sqr{1.5}3}
\def\Z{\hbox{\rm\rlap{Z}\hskip 0.2em Z}}
\def\R{\hbox{\rm I\hskip -0.14em R}}
\def\th #1 #2. #3\par{\medbreak{\bf#1 #2.
\enspace}{\sl#3}\par\medbreak}
\documentclass{llncs1}
\usepackage{setspace}
\usepackage{makeidx}  
\usepackage{graphicx}
\begin{document}
\frontmatter          
\pagestyle{headings}  

\title{On the existence of exotic homotopy 3-spheres}
\titlerunning{Existence of exotic homotopy 3-spheres}  
%
\author{Jorma Jormakka}
\authorrunning{Jorma Jormakka}   
%
\tocauthor{Jorma Jormakka}

\institute{Karhekuja 4, 01660 Vantaa, Finland\\
\email{jorma.o.jormakka@kolumbus.fi}}

\maketitle              

\begin{abstract}
This paper gives a geometric topological proof that exotic homotopy 3-spheres
do not exist.   
\end{abstract}

\vskip 1em
\begin {keywordname}
Homotopy spheres, 3-manifold topology, Morse theory.
\end{keywordname}

\section{Introduction}
A manifold is closed if it is compact and has no 
boundary. In an orientable 3-manifold every 2-dimensional 
embedded submanifold has two sides. Embedding $g$ of a submanifold $N$ 
to a manifold $M$ means that there is a homeomorphism between $N$ and 
$g(N)\subset M$. If $g$ mapping a submanifold $N$ to $M$ is not an
embedding then $g(N)$ has self-intersections. 

An embedded circle $S^1$ is called a loop. If it has self-intersections 
it is only a closed curve. If $l$ and $l'$ are closed curves,
a homomorphism is a continuous mapping $h:I\times I\to M$ such that 
$l=\{h(0,x)|x\in I\}$ and $l'=\{h(1,x)|x\in I\}$. We say that $l$ is 
homotopic to $l'$ and write it $l\simeq l'$. Homotopic images of $S^1$ 
belong 
to the same homotopy class and the set of homotopy classes is a group, 
the fundamental group $\pi_1(M)$ of the manifold $M$. As a special case, if 
$l'$ is a point, $l$ is said to be contractible.
If every loop in the manifold is contractible, $\pi_1(M)=1$ and $M$ is 
said to be simply-connected. 

The Poincar\'e Conjecture states that every
simply-connected closed 3-manifold is homeomorphic to the 
3-dimensional sphere $S^3$.

A Morse function is a function $f:M\to \R$ such that in all
but isolated points there is a diffeomorphism $g$ from a neighborhood $V$
of a point $p\in M$ to $\R^3$ so, that if $g(p)=(x,y,z)$, $f(p)=z+c$
where $c$ is a constant. 

In a closed 3-manifold there are finitely many points where there is no such
homeomorphism. The points are called critical points of $f$ and their 
structure is well known: they are classified by the index $i(p)$. Critical 
points of index 0 are points where there is a homeomorphism $g:V\to \R^3$
from a local neighborhood of a critical point $p$ such that if 
$g(p)=(r,\theta,\phi)$, then $f(p)=r+c$ where $c$ is a constant. A critical
point of index 3 for $f$ is a critical point of index 0 for $-f$, so the
Morse levels $f^{-1}(x)$ are spheres and the level grows towards the center.
A critical point of index 1 is a point $p\in M$ where $f$ has a saddlepoint.
There is a local homeomorphism $g:V\to \R^3$ from a neighborhood $V$ of $p$
which maps the Morse levels $f^{-1}(x)$ to surfaces shown in Figure 1. 
The level $x$ grows to the direction of the arrows in Figure 1. A
critical point of index 2 for $f$ is a critical point of index 1
for $-f$.

\includegraphics {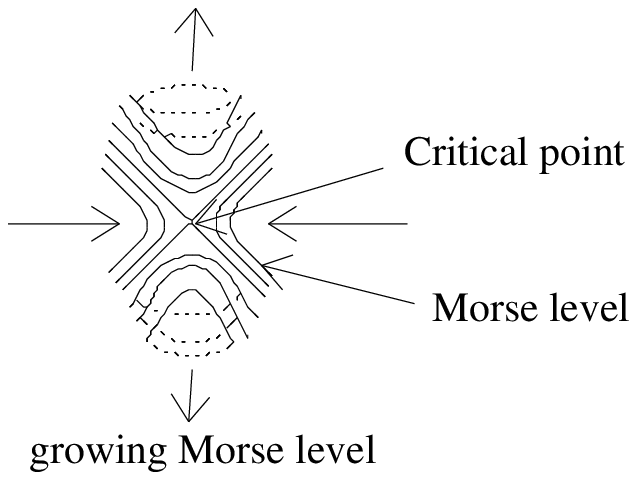}
\hskip 1in	Figure 1.

We will assume that $M$ is a closed and differentiable 3-manifold.

The following notations will be used: $I$ = a closed interval; $S^n$ = a closed n-sphere; $D^2$ = a closed disc;
$T^2$ = a torus; $B^3$ = a closed 3-ball. $D^2_r=\{(x,y)|x^2+y^2\leq r\}$ 
= a closed disc of radius $r$. 
$U^3=\{(x,y,z)|x^2+y^2\leq 1, z\in [0,1], (x,y,1)=(x,y,0)\}$ = a filled torus.
A loop is an image of an embedding $g:S^1\to M$. A path is an image of an 
embedding $g:I\to M$. A disc is an image of an embedding $g:D^2\to M$. 
A collar is an image  of an embedding $g:D^2_1\backslash D^2_{1\over 2}\to M$.

A closed tubular neighborhood $V$ of a loop $l$ is an image of an embedding
$h:U^3\to M$ such that $l$ is the image of the points $x=y=0$ in $U^3$.
A closed tubular neighborhood $V$ of a path $l$ is an image of a homeomorphism
$h:D^2\times I\to M$ such that $l$ is the image of the points $x=y=0$ in $D^2$.
A closed cylinder neighborhood of a disc $D$ is an embedding of $D^2\times I$ to $M$ 
such that $D$ is the image of the points $(r,\theta,1/2)$.
When necessary, we assume that these embeddings are smooth. 

Let $X\approx Y$ mean that the manifold $X$ (possibly with a 
boundary) is homeomorphic to $Y$. Usually we will assume, that 
the homeomorphism is a diffeomorphism. Let $\gamma \simeq \gamma'$ mean that 
the loop $\gamma$ is homotopic to $\gamma'$.  

We will assume, that
the Morse functions are smooth at all points except for a finite set of
critical points. Smooth means here that $f$ is sufficiently many times (say, 3
to be sure) continuously differentiable in local coordinates. Let us write
$M^a_f=\{ x\in M | f(x)\le a \}$, $\partial M^a_f=\{ x\in M | f(x)=a \}$
and if no confusion can arise we write $M^a=M^a_f$. $\partial M^a$ is here 
called a Morse level, but when there is no confusion we also call $a$ a level. 
We will write $M_a=M\backslash int(M^a)$. Let $a_{max}$ and $a_{min}$
denote the highest and lowest levels $a$ respectively such that $f^{-1}(a)$
is not empty.

\section{The idea of the proof}

It suffices to show that every differentiable simply-connected closed 
3-manifold 
is homeomorphic to the 3-dimensional sphere $S^3$ as proving the claim for 
differentiable 3-manifolds proves it for all. As $M$ is simply-connected 
it is orientable, so we start by proving some lemmas for closed, orientable 
and differentiable 3-manifolds.

We do not change the manifold $M$ in the proof but construct different Morse
functions on it. First we recall from a theorem of Smale that it is 
possible to find a Morse
function which has only one critical point of index 0 and of index 3.
In Lemma 2 we show, that we can change the Morse function in a closed tubular 
neighborhood
of a path so, that in the new Morse function critical point of index 2
are on a higher level. We say, that a 2-handle is moved up by changing the
Morse function. In a similar way we can construct a new Morse function 
where a selected critical point of index 1 is on a lower level. We say, 
that a 1-handle is moved down by changing the Morse function.
The modification to the Morse function for moving a 2-handle up or a 1-handle 
down requires changes only in a closed tubular neighborhood of a path.
The modification changes Morse levels in the way that is illustrated 
in Figure 2 by showing a 2-dimensional section. The modification, given 
precisely in Lemma 2, is radially symmetric with respect to the y-axis in Figure 2.

\includegraphics {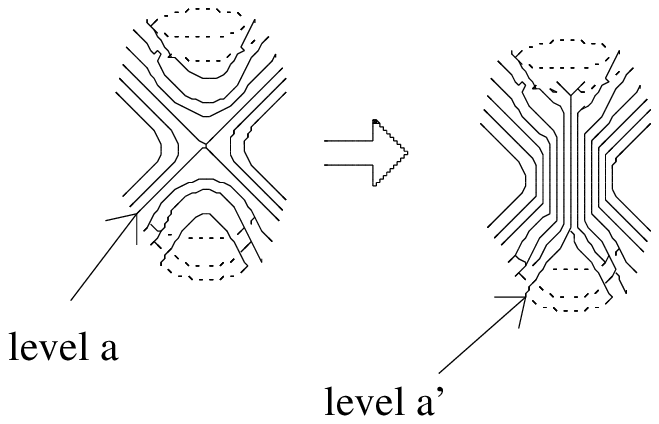}

\hskip 1in	Figure 2.

By moving all 2-handles up and all 1-handles down we get to the situation where
every critical point of index one is on a lower level than any critical point
of index two. 
Let the lowest level of a critical point $p$ of index 2 be $c$ and
the highest level of a critical point of index 1 be $b$.

Let $a$ satisfy $b<a<c$. Then $a$ is a regular level, that is, not a level 
where there are critical points. 
We show that $M^a=\{q\in M|f(q)\le a\}$ and 
$M_a=\{q\in M|f(q)\ge a\}$ are both handlebodies. That is, 3-manifolds with
a boundary obtained by attaching handles to a 3-ball. This is not automatically
true, we must show that the manifolds do not contain embedded
exotic homotopy balls, i.e., 3-manifolds, which 
have the boundary $S^2$ and which are simply-connected but not homeomorphic 
to 3-balls. 

In Lemma 9 we show that inside $M^a$ and $M_a$ there are no exotic homotopy 
balls. This is so 
because if there is an exotic homotopy ball after the addition of 
a 1-handle, there already was an exotic homotopy ball before the addition
of the 1-handle. 
Figure 3 illustrates this: If the blackened part is a part
of an exotic ball, then there already existed an exotic ball before adding the
1-handle because we can connect the parts also by a tube that does not go
through the 1-handle.    

\includegraphics [width=90mm] {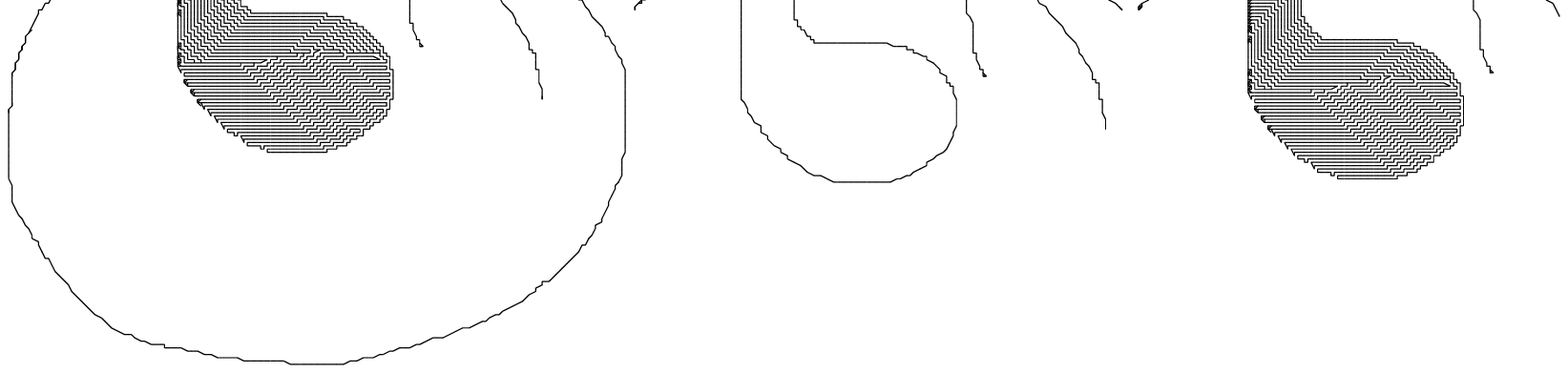}

\hskip 1in	Figure 3.

Let us notice that this is only true for the level $a$ such 
that all 2-handles are above the level $a$ and all 1-handles are below the 
level $a$. On other levels we cannot conclude that there are no 
embedded exotic balls. We cannot even conclude that both sides are
handlebodies with possible embedded exotic homotopy balls, since they may have 
nontrivial fundamental groups. 
We can show these things only for the level $a$. 
In order to show that $M^a$ does not contain exotic balls we use induction 
by adding 1-handles one by one and conclude that no exotic ball exists.
We can show the same for $M_a$ by making a similar induction with the Morse 
function $-f$. 

A Heegaard split is a general way of expressing 3-manifolds by glueing
two handlebodies at their boundaries. The genus $g$ of the split is the
number of handles in the handlebodies. 
In Lemma 10 we show that we have a special kind of Heegaard 
split $M^a$, $M_a$. 
There are homeomorphisms
$g_1$ and $g_2$ such that $g_1(M^a)$ and $g_2(M_a)$ are handlebodies
embedded in $\R^3$, and a boundary homeomorphism $\psi: g_1(M^a) \to g_2(M_a)$.
We can take $g_1$ and $g_2$ differentiable if needed.
Standard noncontractible generators $x_{a,1},...,x_{a,g}$
and contractible generators $y_{a,1},...,y_{a,g}$ can be defined for 
$\partial g_1(M^a)=\partial H_{a,g}$. 
Similarily, the loops $x_{b,1},...,x_{b,g}$, $y_{b,1},...,y_{b,g}$ are 
standard generators for $\partial g_2(M_a)=\partial H_{b,g}$. Figure 4 shows 
these generator loops. 

\includegraphics {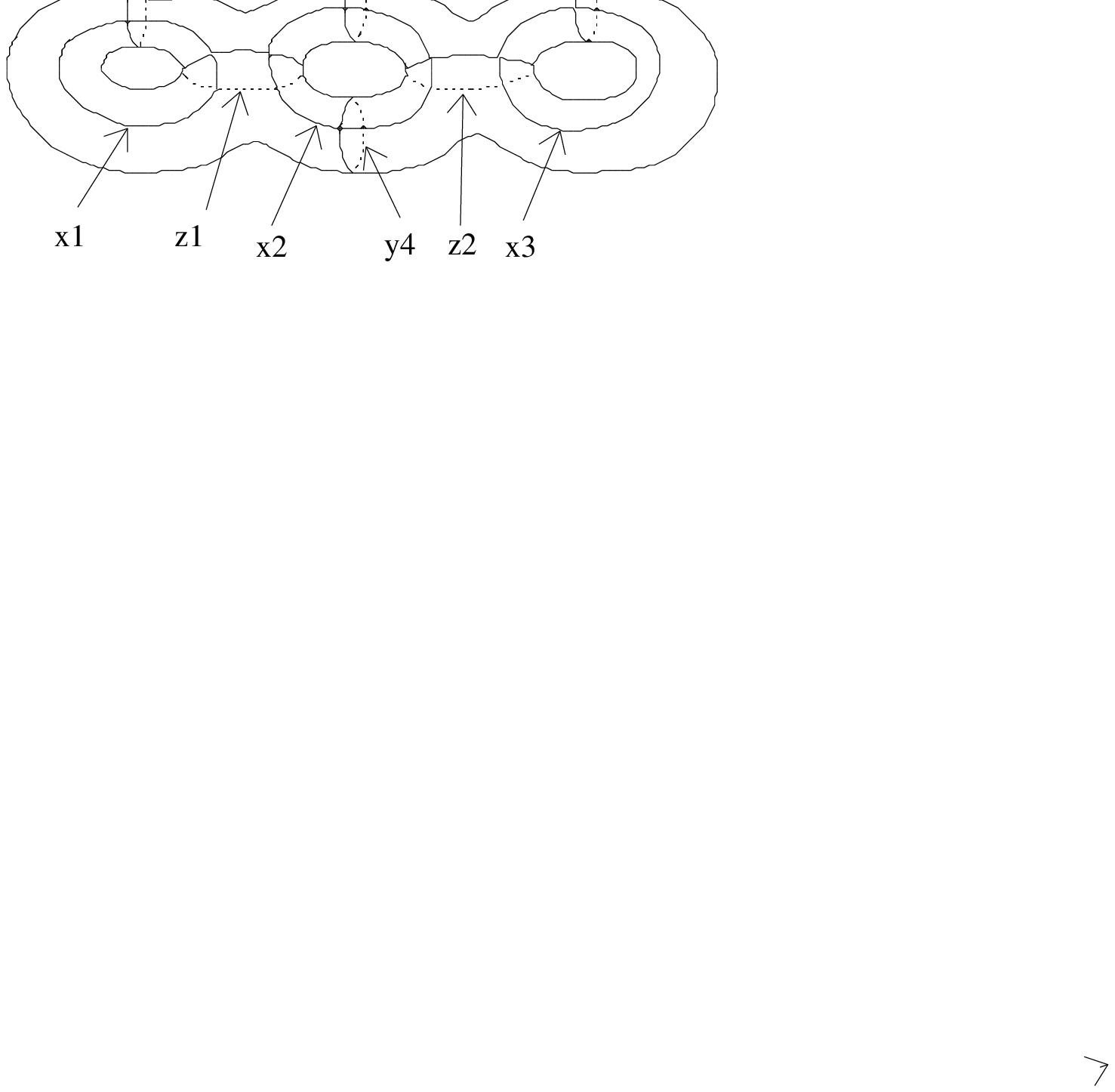}

\hskip 1in	Figure 4.

The boundary homeomorphism $\psi: \partial g_1(M^a) \to \partial g_2(M_a)$ 
induced by a Morse function has the property, that $\psi$ and $\psi^{-1}$ map 
loops to loops.  

The following example for the genus $g=1$ 
shows another important thing. Not all combinations
of standard generators can be represented as loops but only as closed curves 
that have self-intersections. In particular, if a contractible 
generator wraps several times around, there is needed a noncontractible
generator to provide a shift, or the curve is not a loop but has 
self-intersections.  
In Figure 5 a) the surface of a torus $\partial H_{b,1}$
is divided into four squares with side length 4. 
There are identications 
in the horizontal boundaries $((x,4)=(x,-4)$ and in the vertical boundaries
$(y,-4)=(y,4)$ making it a 2-dimensional torus.

\includegraphics[width=100mm]{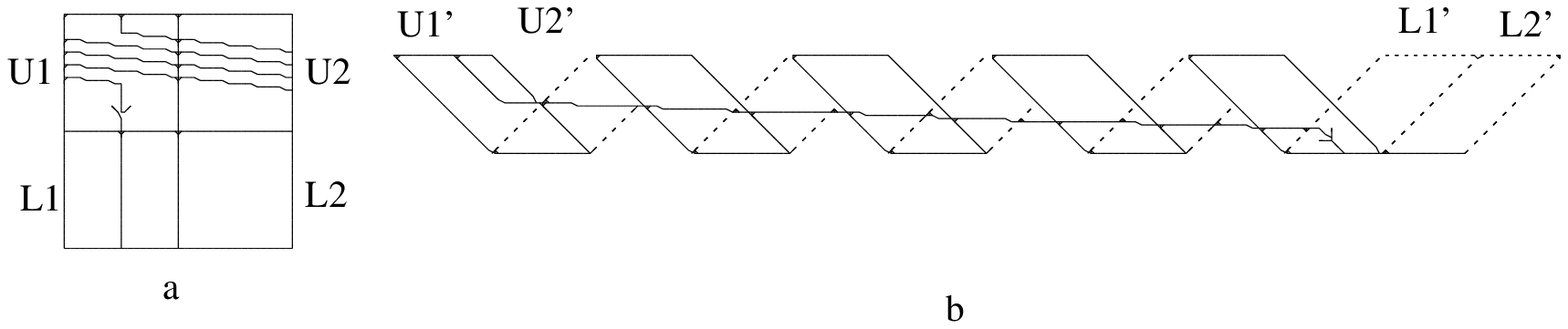}

\hskip 1in	Figure 5.

Let us stretch the part above $y=0$ into a strip of length 16 and width 2,
stretch the part below $y=0$ into a strip of length 4 and width 2 and wrap
the strip 5 times around a torus as in Figure 5 b). A loop drawn in Figure 5
a) maps to a loop drawn in Figure 5 b). Notice, that we need the shift in
Figure 5 a). 

We can define the standard generators on the surface of the handlebody 
$H_{a,1}$ as follows:
$y_{a,1}$ is the loop from $(1,0)$ to $(1,2)$, and $x_{a,1}$ is the loop from
$(1,1)$ to $(21,1)$ in Figure 5 b), and we can define the standard 
generators on the surface of the other handlebody $H_{b,1}$ as follows:
$y_{b,1}$ is the loop from $(2,-4)$ to $(2,4)$, and $x_{b,1}$ is the loop from
$(0,2)$ to $(8,2)$ in Figure 5 a).
The loop drawn to Figure 5 a) is 
$$l_b= (x_{b,1} {y_{b,1}}^{-\alpha} )^4 {y_{b,1}}^{-1+4\alpha}$$
where $\alpha =0.5$ is the offset the loop advances in each round, the 
necessary shift needed to make a loop.

Instead of fundamental groups of $H_{a,1}$ and $H_{b,1}$ we will consider  
the first homology groups $H_1(\partial H_{a,1})$ and $H_1(\partial H_{b,1})$
with integer coefficients. 
The boundary map $\psi:\partial H_{a,1}\to \partial H_{b,1}$ induces a linear
map between the homology groups.  
The first homology group is abelian, so the loop can be expressed 
more simply as
$$l_b=(x_{b,1})^4{y_{b,1}}^{-1}.$$ 
As the group is abelian, we can also use the additive notation:
$$l_b=4x_{b,1}-y_{b,1}.$$

Clearly, we can use linear algebra for solving the generators $x_{a,1}$
and $x_{b,1}$, and if the manifold is simply-connected, we will manage
to solve $x_{a,1}$ and $x_{b,1}$ as linear combinations of $y_{a,1}$ and
$y_{b,1}$. However, this is not only linear algebra with rational coefficients
but there are the 
properties that loops map to loops and not every closed curve can be modified
to be a loop. With these two properties we manage to show in 
Lemma 13 that a loop
$g_1^{-1}(x_{a,1})$ can be associated with a loop $l$ which
is homotopic 
to it in $g_1(M^a)$ and $l$ bounds an embedded disc in $g_2(M_a)$. 
There is a critical point of index 2 in this disc, which we will move down. 

What we want to do with this result is to reduce the genus of the Heegaard 
split by canceling one 2-handle and 1-handle pair. Then induction on the genus
yields the result. The tricky thing is that in order to cancel a pair
we should move a 2-handle down. The modification that is needed for the Morse
function in order to move a 2-handle down is the inverse modification to the
one in Figure 2. The difficult part is that in order to move a 2-handle down
(or 1-handles up, which is the same thing for $-f$) 
we need a closed cylinder neighborhood of a disc. 
In order to move a 2-handle
up (or 1-handles down, which is the same thing for $-f$) 
we only needed a closed tubular neighborhood of a path. Topologically these are
the same and that is why the modification is the same, but while it is easy 
to find a 1-dimensional path which does not meet any critical points, it is not
easy to find a large disc which is avoiding critical points. It is where we
need the simply-connectedness of $M$.  

We can of course find an embedded disc from $M_a$ that bounds any contractible
loop on $\partial M_a$, but this loop has to be moved down to a low Morse 
level. We have to find an embedded collar that has the same boundary as 
this embedded disc. 
If we get such a collar then we can construct a Morse function where
the critical point of index 2 is on a level just above the level of the lowest
critical point of index 1.

We can try to push down a loop that is contractible in $M_a$
and see what happens.
Let us take the loop $l_a=\psi(g_2^{-1}(y_{b,1}))$ on $\partial M^a$
and try to push it down by a family of loops $l_d$
parametrized by $d$ where $d=a$ at the beginning,  

We can move the loop $l_d$ towards the gradient of 
the Morse function $-f$ so that the Morse level of $f$ decreases.
Thus, at each level $d$ holds $l_d\subset \partial M^d$.
In the beginning the loop $l_d$ creates an embedded collar when $d$ decreases
but if $l_d$ meets a critical point of index 1 there may happen several
things - the loop may get knotted, get self-intersections, 
split, merge and so on. 
This creates a complicated situation, and we try to avoid hitting any 
critical points if possible. 

In order to avoid critical points, we try to find a loop
that goes around one 1-handle only in $M^a$, 
so that the loop can be pushed down
creating a collar that does not go through any other critical points. Thus,
we need a loop on $H_{a,1}$ that has only one noncontractible generator 
$x_{a,1}$. If the genus $g=1$ this is of course true, but if $g>1$ we
need to find a suitable loop in Lemma 13. By using linear algebra and
the two properties that loops map to loops and not all closed curves can
be represented as loops we find a suitable loop that has only one 
noncontractible generator $x_{a,1}$. Then we can push the collar down 
all the way without hitting critical points because meeting these
critical points means going around a 1-handle that is created in a critical 
point. 
Our collar goes around the 1-handle that is created on the first 
critical point of index 1. We will not continue the collar to this level
but stop just before the level on some regular level. Thus, our collar avoids
all critical points and stays as an embedded collar. 
 
The embedded collar formed by moving the loop $l_d$ down can be connected
the the disc in order to form a larger disc. We take a closed cylinder 
neighborhood of the disc. Inside the neighborhood
we can the change of the Morse function, using the inverse of the change 
in Figure 2, and move the 2-handle to a lower level. 
The change does not affect the Morse function outside the neighborhood
and will not change the levels of other critical points. 
In this way we can move the 2-handle to a so
low level that the two first critical points after the
0-handle are a 1-handle and a 2-handle and the 2-handle is a disc contracting
the 1-handle. Together they make a 3-ball. We can replace the Morse
function by such a Morse function where these two critical points do
not appear at all. In Lemma 14 we show that this results to a reduction
of the genus of the Heegaard split, which proves by induction
that there are no exotic homotopy 3-spheres.

\section{The proof}

A surface $N\subset M$ is a connected, closed, embedded 2-manifold. For any
regular value $a$ the set $\partial M^a_f$ is a finite set of disjoint 
surfaces. Let $N_a\subset \partial M^a$ be a surface.

\th Lemma 1. Let $M$ be a closed, orientable and differentiable 3-mani\-fold.
Each surface $N_a\subset \partial M^a$ for a regular
value $a$ is 2-sided and orientable, and separates $M$ into two orientable
3-manifolds with boundary $M_1$ and $M_2$, $\partial M_1=\partial M_2=M_1\cap M_2
=N_a$.

\proof As $M$ is orientable the surface $N_a$ separates $M$.$\eop $

\th Theorem D (Smale). Let $M$ be a closed simply-connected differentiable 
3-manifold. 
There exists a Morse function $f:M\to \R$ such that 
$f$ has one critical point of index 0 only and one critical point of index 3.
The number of critical points is finite.

\proof
See the Theorem D of Smale in [2]. $\eop$

\th Lemma 2. Let $M$ be a closed orientable and differentiable 3-manifold. 
Let $f:M\to \R$ be a Morse function such that $f$ has the critical points 
$p_0,...,p_n$.

{\sl Let $k$ be such that $i(p_k)=2$. Let ${a'}_1>a_1=f(p_k)$,
${a'}_1<a_{max}$, be a noncritical 
level such that there exists a path $l:I\to M$ satisfying }

{\sl $l(0.5)=p_k$, }

{\sl $f(l(1))>{a'}_1$, }

{\sl for each $x,y\in (0,1)$ holds that if $x<y$ then $f(l(x))<f(l(y))$ and}

{\sl $l(I)\cap \{ p_m | 0\le m\le n\}=\{p_k\}$.}

{\sl Then there is a Morse function $f':M\to \R$ such that:}

{\sl $f'$ and $f$ have the same critical points and they have the same indices,}

{\sl $f'(p_k)={a'}_1$ and $f'(p_m)=f(p_m)$ for all critical points $m \not =k$,}

\proof
Select a small closed tubular neighborhood 
$V$ of $l(I)$ such that $V\cap \{p_m | 0 \le m \le n\}=\{p_k \}$. 

Define 
$A=\{ (r,\theta,z)\in \R^3 | 0 \le r \le 1, 0\le \theta < 2\pi , 0 \le z \le 1\}$.
Select $s\in (0,1/2)$. Let ${w'}_s:A\to A$ be defined as 

${w'}_s(r,\theta,z) = (r,\theta,z+(1-r)(0.5-s){z \over s})$ 
if $0\le z\le s$, $0\le r\le 1$,

${w'}_s(r,\theta,z) = (r,\theta,z+(1-r)(0.5-s))$ 
if $s< z\le 0.5$, $1-{{0.5-z}\over {0.5-s}} \le r \le 1$.

We extend ${w'}_s$ to $1\ge z>0.5$ by:

${w'}_s(r,\theta,z) = (r,\theta,1-w_s(r,\theta,1-z))$ 
if $0\le 1-z\le s$, $0\le r\le 1$,

${w'}_s(r,\theta,z) = (r,\theta,1-w_s(r,\theta,1-z))$ 
if $s< 1-z\le 0.5$, $1-{{0.5-(1-z)}\over {0.5-s}} \le r \le 1$.

Let $w_s:A\to A$ be a smooth fuction for $r<1$ approximating ${w'}_s$. 
Let us find a diffeomorphisms $g$ of $A$ to the tubular neighborhood $V$ so that

the points with $z=0.5$ map to $l(I)$, $g(0,0,0.5)=p_k$, 
 
$g(0,0,0)=l(0)$, $g(0,0,1)=l(1)$ and

$g^{-1}(\partial {M^a}_f\cap V)=\{ w_{a_1}(r,\theta ,a)|0\le r\le 1,0\le \theta<2\pi \}$.

Let us write $A_{a,a_1}=g^{-1}(\partial {M^a}_f\cap V)$.

Let us construct $f':M \to \R$ as follows:

$f'|_{M\backslash V}=f|_{M\backslash V}$,

$f':V\to \R$ is defined by the condition: if $p\in g^{-1}(A_{a,{a'}_1})$ 
then $f'(p)=a$.

Then $f'(p_k)={a'}_1$ since $g^{-1}(p_k)=(0,0,0.5)\in A_{{a'}_1,{a'}_1}$.

Since ${a'}_1>a_1$, $f'$ does not have other critical points than $p_k$ in
$V$. (If ${a'}_1<a_1$, $f'$ could have critical points in $\partial V$ since 
there could be created a new saddlepoint.) Finally smoothen $f'$ on
$\partial V$. $\eop$

Figure 2 shows the construction of $f'$ in $V$. We simply wrote down an
explicit function for doing the modification in Figure 2.

\th Lemma 3. Let $M$ be a closed simply-connected differentiable 
3-manifold. There exists a Morse function $f:M\to \R$ such that $f$ has
one critical point of index 0 only, one critical point of index 3 only and
on each critical level there is one critical point only.

\proof
By Theorem D of [2] there is a Morse function $f$ such
that $f$ has one critical point of index 0 only and one critical point of
index 3 only.

By applying the modification in Lemma 2 successively we can construct a 
Morse function which has the same critical points (and of same indices) 
as $f$ but such that all critical points are on distinct levels. $\eop$

Let us denote the minimum and maximum values of $f$ by 
$f(p_0)=a_{min}$ and $p_n=a_{max}$. 

\th Lemma 4. Let $M$ be a closed simply-connected differentiable 
3-manifold. If a Morse function $f:M\to \R$ has
one critical point of index 0 only and one critical point of index 3 only, 
then for each level $a$ $\partial M^a$ is connected.

\proof If $\partial M^a$ is not connected, then there exists a noncontractible loop
passing the critical point of index 0, the critical point of index 3 and two
points in two components of $\partial M^a$. $\eop$

\th Lemma 5. Let $M$ be a closed simply-connected and differentiable 
3-manifold. 
Let $f:M\to \R$ be a Morse function such that $f$ has the critical points 
$p_0,...,p_n$.

{\sl Let $k$ be such that $i(p_k)=2$. Let ${a'}_1>a_1=f(p_k)$ be a noncritical 
level such that ${a'}_1<f(p_n)=a_{max}$. }

{\sl Then there is a Morse function $f':M\to\R$ such that:}

{\sl $f'$ and $f$ have the same critical points and they have the same indices,}

{\sl $f'(p_k)={a'}_1$ and $f'(p_m)=f(p_m)$ for all critical points $m \not =k$,}

\proof

If there exists a path $l:I\to M$ as in Lemma 2, then the result follows 
from Lemma 2. The only problem is to select $l:I\to M$ such that

$l(I)=\{ p_m | 0\le m\le n\}=\{p_k\}$

and that the path satisfies:
for each $x,y\in (0,1)$ holds that if $x<y$ then $f(l(x))<f(l(y))$.

The points $p_0$ and $p_n$ are on levels outside $l(I)$. If $i(p_m)=1$, 
then $p_m$ can be avoided as a critical point of index 1 in the direction 
of growing Morse levels does not force a path to go through the critical point, 
see Figure 1. Let $l(I)$ be selected so that it does not
intersect with critical points of index 1.

It is not obvious that a critical point of index 2 can be avoided when moving 
a 2-handle in the direction of growing Morse level. Let $p_m$ be 
a critical point of index 2 with $f(p_m)>f(p_k)$, $p_m\in l(I)$ and
there is no critical point $p_j\in l(I)$ such that $f(p_k)<f(p_j)<f(p_m)$.

Assume that there is no path which omits $p_m$. Select a small
embedded circle in a neighborhood of $p_k$ on the level $f(p_k)+\epsilon$ 
such that a path from every point of the circle a path towards grad($-f$) 
passes through $p_k$. Continuing the circle to the direction of $grad(f)$
defines a surface $N\subset M$ such that $N=g(I\times I)$ where 

$g(0,x)=p_k$, $g(1,x)=p_m$ for $x\in [0,1]$,
 
if $y>z$, $y,z\in [0,1]$, then $f(g(y,x)) > f(g(z,x))$ for any $x\in [0,1]$,

$g(y,0)=g(y,1)$ for any $y\in [0,1]$.

If there is no such $N$, then $l(I)$ can be selected to omit $p_m$.

$N$ is embedded as in all points $x\in N$, $x\not\in \{p_m,p_k\}$
the gradient of $f$ is unquely defined. 

The surface $N$ separates $M$ as $M$ is orientable. As $p_m$ is a saddlepoint,
there are two points $p$ and $q$ on the opposite sides of $N$ such that both
$p$ and $q$ belong to $\partial M^{f(p_m)+\epsilon}$ for some small 
$\epsilon>0$. Let $a=f(p_m)+ \epsilon$. By Lemma 4 $\partial M^a$ is 
connected, therefore there is a path from $p$ to $q$ in $\partial M^a$. 
This path must intersect $N$ as $N$ separates $M$. However, all points 
$x\in N$ have $f(x)\le f(p_m)<a$. This contradiction shows that there is a 
path $l(I)$ omitting $p_m$. $\eop$

\th Lemma 6. Let $M$ be a closed simply-connected and differentiable 
3-manifold. 
Let $f:M\to \R$ be a Morse function such that $f$ has the critical points 
$p_0,...,p_n$.

{\sl Let $k$ be such that $i(p_k)=1$. Let ${a'}_1<a_1=f(p_k)$ be a noncritical 
level such that ${a'}_1>f(p_0)=a_{min}$. }

{\sl Then there is a Morse function $f':M\to\R$ such that:}

{\sl $f'$ and $f$ have the same critical points and they have the same indices,}

{\sl $f'(p_k)={a'}_1$ and $f'(p_m)=f(p_m)$ for all critical points $m \not =k$,}

\proof 
A critical point of index 1 of $f$ is a critical point of index 2 of $-f$.
The result follows from Lemma 5. $\eop$

\th Lemma 7. Let $M$ be a closed 
simply-connected and differentiable 3-manifold. 
There exists a Morse function $f:M\to \R$ such that

{\sl $f$ has the critical points $p_0,...,p_n$,}

{\sl the indices of the critical points are $i(p_0)=0$, $i(p_n)=3$,
$i(p_m)=1$ for $0<m\le l$, $i(p_m)=2$ for $l<m<n$,}

{\sl for any $m<j$ holds $f(p_j)<f(p_m)$.}

\proof 
By Lemma 5 there is a Morse function $f'$ such that one critical point $p_k$
of index 2 is on a regular level $a$ with $a_{max}>a>f(p_m)$ for $m<n$, 
$m\not = k$.
  
By applying Lemma 5 sequentially to all critical points of index 2
we find after $n-l-1$ steps a Morse function where critical points of index 2
are on higher levels than all critical points of index 1 or 0. $\eop$

\th Lemma 8. Let $M$ be a closed simply-connected and differentiable 
3-manifold. 
Let $f:M\to \R$ be a Morse function such that

{\sl $f$ has the critical points $p_0,...,p_n$,}

{\sl the indices of the critical points are $i(p_0)=0$, $i(p_n)=3$,
$i(p_m)=1$ for $0<m\le l$, $i(p_m)=2$ for $l<m<n$,}

{\sl for any $m<j$ holds $f(p_j)<f(p_m)$.}

{\sl Let $a$ be a regular level. If $a<f(p_{i+1})$, then $\pi_1(M^a)$
is free. If $a>f(p_i)$, then $\pi_i(M_a)$ is free.}

\proof
By the Morse theory $\pi_i(M^a)=1$ for $a\in (a_{min}, f(p_1))$. In each
critical point $p_k$ of index 1 ($0<k\le i$) by the Morse theory
$\pi_1(M^{f(p_k)+\epsilon})=\Z \times \pi_1(M^{f(p_k)-\epsilon})$ for
a small $\epsilon>0$. By induction $\pi_1(M^a)$ is free for any
regular level $a<f(p_{i+1})$. As critical points of index 2 for $f$
are critical points of index 1 for $-f$, $\pi_1(M_a)$ is free
for any regular level $a>f(p_i)$. $\eop$

A 1-handle is a cylinder $H\approx D^2\times I$. A handlebody of genus $k$ 
is by definition obtained by adding $k$ 1-handles to $B^3$ so, that the 
1-handle $H$ is fixed at the discs $D^2\times \{0,1\}$ to the boundary
of $B^3$. At a critical point of index 1 a 1-handle 
is added to $M^a$. An exotic homotopy ball is an embedded
manifold $B^3_e$ which has the boundary $\partial B^3_e\approx S^2$ 
and $\pi_1(B^3_e)=1$.

\th Lemma 9. Assume that a 3-manifold $M_1$ (with boundary) is obtained
from a 3-manifold $M_2$ by addition of a 1-handle $H\approx D^2\times I$,
$M_1=M_2\cup H$. Assume, that $M_1$ and $M_2$ are orientable. If $M_1$ 
contains an exotic homotopy ball $B^3_e$. Then $M_2$ contains an embedded
exotic homotopy ball  ${B'}^3_e$.
 
\proof
The embedding of $B^3_e$ 
to $M_1$ is not wild as both manifolds $M_1$ and $B^3_e$ have the same 
dimension 3. Wild embeddings happen only with submanifolds of lower dimension.
As $B^3_e$ is not wildly embedded, $\partial B^3_e\approx S^2$ is not
a wildly embedded sphere.

Let $F$ be a smoothly embedded disc cutting the 1-handle $H$. $F$ can be 
assumed to be in a general position so, that $F\cap \partial B^3_e$ is
a set of disjoint circles $A_0=\{C_i | i\in I_0\}$, $C_i\approx S^1$, $i_0$ finite.

The set $A_0$ contains a nonempty subset $B_0=\{S_j |j\in I_1\}\subset A_0\}$
where each $S_j\subset \partial B^3_e$ 
is a boundary of a disc $D_j \subset \partial B^3_e$ satisfying
$int (D_j) \cap A_0\}=\emptyset$. In order to see that $B_0$ is nonempty,
notice, that each $C_i$ separates $\partial B^3_e\approx S^2$. There must be
such circles on $\partial B^3_e$ that one of the separated sides does not 
contain smaller circles. Such a circle $C_i$ is a circle $S_j$.

Take one $S_j$ and a circle ${S'}_j\subset F$ which is slightly bigger than 
$S_j$. There is a disc ${D'}_j\subset F$ such that $\partial {D'}_j={S'}_j$.
There may be other circles $C_i$ than $S_j$ inside ${D'}_j$ but that does 
not matter here. Replace the disc ${D'}_j$ by a disc ${D''}_j$ which is close 
to $D_j$. Then $D_j$, ${D''}_j$ and the annulus $F_0\subset F$ between the
circles $S_j$ and ${S'}_j$ separates $M_1$ and one side is a 3-ball $B_0$, it is
the side which does not contain points of $int (B^3_e)$. 
The disc ${D''}_j$ separates $M_2$ into two components. 
Let $U$ be the component which contain points of $int (B_0)$.

Let $V$ be a closed neighborhood of $F$ in $H$ and $V_1$ a smaller closed
neighborhood of $F$ in $V$. Then $F$ separates $V$ and $V_1$ and $\partial V$ 
is the union of a collar, which is subset of $\partial H$, and two discs
$D_a$, $D_b$ which have boundary at $\partial H$. 

By construction ${D''}_j$ and $D_j$ insersect with $F$ only at 
${S'}_j$ and $S_j$. Therefore we can select $V$ so small, that 
the component in $V\cap B^3_e$ which contains $D_j$ is a collar
$C=c(S^1\times I)$, for an embedding $c:S^1\times I \to H$ such that
$C\subset V$, $c(S^1\times \{0\})\subset D_a$, $c(S^1\times \{1\})\subset D_b$,
$C\cap F=S_j$, $C\subset \partial B^3_e$.
 
Similarily we can find a collar
$C'=c'(S^1\times I)$, for an embedding $c':S^1\times I \to H$ such that
$C'\subset V$, $c'(S^1\times \{0\})\subset D_a$, $c'(S^1\times \{1\})\subset D_b$,
$C'\cap F={S'}_j$, $C'\cap \partial B^3_e=\emptyset$ and
${D''}_j\cap V\subset C'$. ${D''}_j\cap V$ is on one side of $F$ in $V$, 
let that be the side where $D_a$ is. 
 
Let us define a homeomorphism $g':M_1\to M_1$ so, that $g'$ keeps the 
disc ${D''}_j$ fixed, expands the 3-ball $B_0$ to ${B'}_0$so, that the 
disc $D_j$ is pushed to the component $K_b$ of $V\backslash int (V_1)$ which contains $D_b$.
The mapping $g'$ pushes a part $A$ of the homotopy ball $B^3_e$ 
which was in $M_2$ to $K_b$.

We must move this part of $B^3_e$ back to $M_a$ but not through $F$.
Let $l$ be a path from a point in $V_1\cap K_b$ to a point in $V_1\cap K_a$
which does not go through $F$. We can select $l$ so that it does not
intersect with $B^3_e$ as $\partial B^3_e\approx S^2$. Let $V_2$
be a small closed tubular neighborhood of $l$ in $M_1$. We 
connect the tube $V_2$ to $A$ and move $A$ through the tube $V_2$ to
the other component $K_a$ of $V\backslash (V_1)$ 
by expanding a 3-ball into $K_b$.
Finally we restore $A$ back to $M_2$ by decreasing the 3-ball in $U$.
Now the tube $V_2$ is filled with a 3-ball, $M_2$ is restored and the
handle $H$ is again a handle. Let us call this homeomorphism of $M_1$
to $M_1$ by $g''$. See Figure 3 where the mapping is shown. 
  
Define $g$ as the combined homeomorphism $g:M_1\to M_1$, $g(p)=g''(g'(p))$. 
Then $g(M_2)=M_2$. 

Let us rename the exotic homotopy ball $g'(B^3_e)$ as $B^3_e$.

There are now less circles in the intersection set 
$A_1=F\cup \partial B^3_e$. If $A_1$ is nonempty, 
let us return to the step of selecting one $S_j$, repeat the procedure and
get a smaller set $A_2$.

We repeat the procedure as long as $F\cup \partial B^3_e \not = \emptyset$.
When the set is empty, we have obtained an exotic homotopy ball which
is contained in $M_2\cup K_1\cup K_2\approx M_2$. So, there is an exotic 
homotopy ball ${B'}^3_e$ in $M_2$. 
$\eop$
     
\th Lemma 10. Let $M$ be a closed simply-connected differentiable 
3-manifold. Then there exists a Morse function $f:M\to \R$ such that

{\sl $f$ has the critical points $p_0,...,p_n$,}

{\sl the indices of the critical points are $i(p_0)=0$, $i(p_n)=3$,
$i(p_m)=1$ for $0<m\le l$, $i(p_m)=2$ for $l<m<n$,}

{\sl for any $m<j$ holds $f(p_j)<f(p_m)$.}

{\sl For any regular level $a$ holds: if $a<f(p_{l+1})$, then $M^a$
is a handlebody and if $a>f(p_l)$, then $M_a$ is a handlebody.}

\proof
By Lemma 7 there exists such a Morse function. By Lemma 8 the fundamental
groups of $M^a$ and $M_a$ are free. A compact orientable 3-manifold 
with a nonempty connected boundary, which has a free fundamental group 
is a handlebody with embedded homotopy balls.

If there is an embedded exotic homotopy ball in $M^a$ for a regular 
$a<f(p_{l+1})$, then Lemma 9 used $l$ times shows 
that there is an exotic homotopy ball in $M^a$ for $a_{min}<a<f(p_1)$.
This is a contradiction as for $a_{min}<a<f(p_1)$, $M^a$ is a 3-ball.
This shows that there are no embedded exotic homotopy balls in $M^a$. 
Therefore $M^a$ is a handlebody if $a<f(p_{l+1})$.

The same argument for $-f$ shows that there are no embedded exotic homotopy 
balls in $M_a$ for $a>f(p_l)$. 
Therefore $M_a$ is a handlebody if $a>f(p_l)$. $\eop$  

Let us now summarize the results obtained so far. These results do not require
the manifold $M$ to be simply-connected, and the results are very natural.
We have created a Heegaard split for the manifold, but this Heegaard split
has the advantage that it is induced by a Morse function and we learned how to
move handles by Lemma 2. 

\th Definition 1.
Let $M$ be a closed, orientable and differentiable 3-manifold. 
Let $f$ be a Morse function with critical points $p_0,...,p_n$.
Let $a$ be a regular level with $f(p_m)<a$ for all critical points $p_m$ 
with $i(p_m)\le 1$ and $a<f(p_m)$ for all critical points $p_m$ with $i(p_m)>1$. 
The 3-manifold $M$ has a Heegaard split induced by the Morse function $f$
if there there are homeomorphisms $g_1$, $g_2$ of $M^a$ and $M_a$ to standard
handlebodies $g_1(M^a)$ and $g_2(M_a)$ in $\R^3$
and a boundary identification homeomorphism $\psi$ such that the following diagram 
commutes

$$\vbox{\rm\halign {\hfil #\hfil && \quad \hfil #\hfil \cr
 $M^a$&  ${{g_1} \above 0pt {\to}}$ & $g_1(M^a)$ & $\hookrightarrow$ & $\R^3$\cr
 $\uparrow$  & & $\uparrow$\cr
 $\partial M^a$  & & $\partial g_1(M^a)$\cr
 $id$ $\downarrow$& &$\psi$ $\downarrow$\cr
 $\partial M_a$  & & $\partial g_2(M_a)$\cr
 $\downarrow$ &  & $\downarrow$\cr
 $M_a$&  ${{g_2} \above 0pt {\to}}$ & $g_2(M_a)$ & $\hookrightarrow$ & $\R^3$\cr
}}$$

\th Corollary 1.
Let $M$ be a closed differentiable and
simply-connected 3-manifold. Let $f$ be a Morse function as in Lemma 10.
Let $a$ be a regular level with $f(p_i)<a<f(p_{i+1}$, then
$M^a$, $M_a$ is a Heegaard split induced by the Morse function $f$.

\proof
The claim follows directly from Lemma 10. $\eop$

The homeomorphism $\psi$ of a Heegaard split is called a 
gleuing mapping. It is sufficient within a homeomorphism to define 
the glueing mapping by defining how the standard generators are mapped. 
A Heegaard split induced by a Morse function is a Heegaard split but 
our goal is to use the function $f$ and to move a 2-handle down. 

The following property is useful.

\th Lemma 11. Let $\psi$ be the glueing mapping of a Heegaard split induced by
a Morse function $f$. If $l\in \partial M_a$ is a loop, then 
$l'={g_1}^{-1}(\psi^{-1}(g_2(l)))$ is a loop in $\partial M^a$.
If $l\in \partial M^a$ is a loop, then 
$l'={g_2}^{-1}(\psi(g_1(l)))$ is a loop in $\partial M_a$.

\proof
The Morse function $f$ has $a$ as a regular level. Therefore
the glueing mapping $\psi:\partial M^a\to \partial M_a$ is bijective. 
It follows that if $l$ is a loop, $l'$ is a loop. $\eop$

A Heegaard split of genus 1 can be defined 
by glueing two standard handlebodies $H_{a,g}$ and $H_{b,g}$ of genus g
together so, that the standard generators
$x_{a,i},y_{a,i}$ of $H_{a,g}$ are mapped to words in generators
$x_{b,1},y_{b,1}$ of $H_{b,g}$ by $\psi$.
A fundamental group of the boundary $\partial H_g$ of a handlebody $H_g$ 
is a group generated by 
$$x_1,...,x_g,y_1,...,y_g$$
with one relation $x_1y_1x_1^{-1}y_1^{-1}\cdots x_gy_gx_g^{-1}y_g^{-1}$.

\th Lemma 12. Let $M$ be simply-connected and obtained as a Heegaard split of 
genus $g$. The boundary map 
$\psi:\partial H_{a,g}\to \partial H_{b,g}$ 
induces a linear map between the first homology groups 
$H_1(\partial H_{a,g};\Z)$ and $H_1(\partial H_{b,g};\Z)$. 
These groups are generated by $2g$ generator loops  
$$x_{c,1},...,x_{c,g},y_{c,1},...,y_{c,g} \hskip 2em c\in \{a,b\}$$ 
where $x_{c,j}$ are noncontractible in the handlebody $H_{c,g}$
and $y_{c,j}$ are contractible in the handlebody $H_{c,g}$.
Every $x_{c,j}$ can be expressed as
$$x_{c,j}=A_j{Y_a}^T+B_j{Y_b}^T \eqno(2)$$
where $A_j$ and $B_j$ are row vectors with rational entries
and $Y_c=[y_{c,1},\dots, y_{c,g}]^T$. 

\proof
The first homology group is the abelianization of the fundamental group and
thus it is abelian and torsion free. We can use the same generators for
the first homology group, or select another set of generators.  
As the group is abelian, we can use an additive notation and write a
word in the generators as

$w=\sum_i (a_ix_i+b_iy_i)$.

The pullback of the boundary map $\psi:\partial H_{a,g}\to \partial H_{b,g}$ 
is a linear map between the homology groups.  
Thus, the map between the first homology groups 
$H_1(\partial H_{a,g};\Z)$ and $H_1(\partial H_{b,g};\Z)$. 
induced by the glueing mapping 
can be described in a matrix form as

$Z_a=SZ_b$ , where 

$Z_c=[X_c,Y_c]^T=[x_{c,1},...,x_{c,g},y_{c,1},...,y_{c,g}]^T$,
$c\in \{a,b\}$. 

The entries of the matrix $S$ are rational numbers.
There are 2$g$ equations to solve $x_{a,j}$ and $x_{b,j}$.
If the matric $S$ is invertible,
every $x_{c,j}$ can be expressed as

$$x_{c,j}=A_j{Y_a}^T+B_j{Y_b}^T$$

where $A_j$ and $B_j$ are row vectors with rational entries.
Since $y_{c,j}$ is contractible in $H_{c,g}$ ($c\in \{a,b\}$),

$y_{a,j}\simeq 1$ for all $j=1,...,g$ in $H_{a,g}$,

$y_{b,j}\simeq 1$ for all $j=1,...,g$ in $H_{b,g}$.

It follows that if the matrix is invertible, then 

$x_{a,j}\simeq 1$ for all $j=1,...,g$ in $M$,

$x_{b,j}\simeq 1$ for all $j=1,...,g$ in $M$.

This is as we expect since the manifold $M$ is simply-connected.

However, if the matrix $S$ is not invertible, then we cannot solve every
$x_{c,j}$ as a linear combination of $y_{a,k}$ and $y_{b,k}$.
In this case there are too few linearly independent equations and at least
one $x_{a,j}$ can only be expressed as a funtion of  
$y_{a,k}$, $y_{b,k}$, $x_{b,k}$ and $x_{b,m}$ for some indices $k$ and 
$m\not=j$. In this case $x_{a,j}$ cannot be contractible in $M$. Then
$M$ is not simply-connected. We conclude that the matrix $S$ must be
invertible. 

In order to illustrate this conclusion by an example, let us take $g=1$ and
let the mapping be $x_{a,1}=x_{b,1}$, $y_{a,1}=y_{b,1}$. Here we cannot 
solve $x_{a,1}$ as a linear combination of $y_{a,1}$ and $y_{b,1}$ and
the manifold is also clearly not simply-connected.   
$\eop$

\th Lemma 13. If $M$ obtained as a Heegaard split of genus $g$
induced by a Morse function is simply-connected and $x_{a,j}$ is
a standard noncontractible generator loop of $H_{a,g}=g_1(M^a)$ 
then there exists 
a loop $l$ homotopic to $x_{a,g}$ in $H_{a,g}$ which bounds a disc
in $H_{b,g}$.

\proof
We will select the generators of the first homology groups of 
$\partial H_{c,g}$, $c\in \{a,b\}$ slightly differently as $y_j$ in Figure 4.
Let us first take a large filled torus
and let its noncontractible standard generator
be called $x_{c,1}$ and contractible standard generator be called $y_{c,1}$.
Then we glue $g-1$ small handles to the outer sphere of this filled torus.
This gives a handlebody of genus $g$. We select $x_{c,j}$ as noncontractible
generators as in Figure 4, but the contractible generators $y_{c,j}$ we
select such that they are on the surface of the big torus. Thus, they go
through the holes created by the small handles. In Figure 4 we have marked
these loops by $z_j$. 
We can use these generators
in Lemma 12 and 
as $M$ is simply-connected, there is a solution    

$x_{a,j}=A_j{Y_a}^T+B_j{Y_b}^T$

where $A_j$ and $B_j$ are row vectors with rational entries. 

In the example of Figure 5 a) we needed a shift $\alpha$ 
in order to create a curve without making self-intersections, that is, to
make a loop. 

Moving the generators of $H_{a,g}$ to the left-hand side and  
the generators of $H_{b,g}$ to the right-hand side gives  

$x_{a,1}-A_1{Y_a}^T=B_1{Y_b}^T$

Let the vector $B_1$ have the components 
$[b_1,\dots, b_g]$. Six observations can be made. 

1) If the coefficients of $A_1$ are integers, the left-hand side can be 
presented as a loop. That is, we can use $x_{a,1}$ to shift multiple
windings of $y_{a,k}$ in such a way that there are no intersections, 
similarly as is done in Figure 5 a).  

2) Assume that the coefficients of $A_1$ are rational numbers but not integers.
In that case the left-hand side is a set of disjoint paths that have no 
self-intersections. We can find an integer $k$ such that when 
the equation is 
multiplied by $k$, the left-hand side can be presented as a loop. 

3) It is generally true that if there is no shift, we do not get a loop 
if a curve winds up multiple times. This can be seen in Figure 5 a). 
Thus, if $b_j$ is an integer, the curve $b_jy_j$ on the surface 
of a handlebody is not a loop unless every nonzero $b_j=\pm 1$. 
It follows that the 
right-hand side can be presented as a loop only if every nonzero $b_j=\pm 1$.  

4) In the case 2) we had to multiply the equation by an integer $k$ in order
that the left-hand side could be presented as a loop. In order to see that
we get a loop in this case proceed as follows. Consider Figure 4 and notice 
that we could add any number of handles to the torus in the center of Figure 4.
Thus, we have the generator $x_2$ and a number of generators $z_j$. 
Take a loop where $x_2$ appears one time only and $z_j$ 
appears $c_j$ times where $c_j\in \Z$. 
Cut the central torus in Fig. 4 along a loop parallel to $x_2$ 
and deform it into a strip $S^1\times [0,1]$. 
Glue two copies of the strip together by identifying $S^1\times \{0\}$
with $S^1\times \{1\}$ of the second copy. Repeat this $k$
times and finally identify the sides $S^1\times \{0\}$ 
and $S^1\times \{1\}$. This yields a loop where $x_2$ 
appears $k$ times and $z_j$ appears $c_j$ times. This loop is mapped to a loop
and because of 3) every nonzero $b_j$ after multiplying by $k$ is $\pm 1$. 
It follows that the integer $k$ is $\pm 1$. 

5) As the left-hand side is not contractible in $H_{a,g}$,
the right-hand side has at least one coefficient $b_j$ which is not zero.

6) Loops map to loops in $\psi$, therefore the right-hand side is a loop
in $g_2(M_a)$. As it
is loop, every $y_{b,i}$ must have coefficient $\pm 1$ as there is
no term $x_{b,i}$. So

$l'=B_1{Y_b}^T=\sum_{1\le i\le g}b_iy_{b,i}$

where $b_i=\pm 1$ or zero. 
This loop is contractible and bounds an embedded 
disc in $H_{b,g}$. 
 
We conclude that the preimage $g_1^{-1}(l)$ is
a loop in $\partial M^a$, homotopic to $g_1^{-1}(x_{a,1})$ in $M_a$
and bounds an embedded disc in $M_a$. $\eop$

\th Lemma 14. Let $M$ and $f$ be as in Lemma 10 and let $M$ have a Heegaard
split of genus $g$ induced by $f$. There is a Morse function $f'$ such that
$f'$ satisfies the conditions of Lemma 10 and the Heegaard split induced by
$f'$ has genus $g-1$.

\proof
Let $x_{a,1}$ be as in Lemma 13. By Lemma 13 there is a Morse fuction
$f$ satisfying Lemma 10 such that 
$x_{a,1}$ bounds a disc in the handlebody $H_{b,g}$. We can assume that
the loop corresponding to the generator $x_{a,1}$ is created by the 
addition of the first 1-handle in a critical point $p_1$. 
On a regular Morse level $c$ where $c$ is only slightly higher than 
$f(p_1)$ the manifold is divided into two parts $M^c$ and $M_c$. The
boundaries of $M^c$ and $M_c$ are toruses. The lower manifold $M^c$ is
a filled torus but we cannot know if the upper manifold $M_c$ is a filled
torus. 
Fortunately, we do not need to
know that. We are only interested in the gleuing mapping between 
these manifolds.
If $M_c$ is not a filled torus, we replace it by a filled torus.
Now we have two filled toruses glued together by a boundary map. Thus, 
we have obtained a Heegaard split of genus 1, though not for the original 
manifold. Every simply connected closed 3-manifold 
which has a Heegaard split of genus 1 is homeomorphic to $S^3$. This means 
that the modified $M_c$ and the original $M^c$ make $S^3$. That means that
the boundary map must be trivial: a noncontractible generator is mapped to a 
contractible generator. 
Let us now restore $M_a$. It was only temporarily replaced by a filled torus
in order to investigate the boundary mapping. The boundary mapping is of course
still the same. We know now that the boundary mapping is trivial. 

We have the disc in $H_{b,g}$ that has $l$ as a boundary,
and a collar from $l$ to the level $c$. This collar does not meet any critical
points. We can make a similar move of a handle as is done in Figure 2
but this time we move a 2-handle down. 

The collar can be added to the disc to make a bigger disc.
The disc contains a critical point $p_m$ of index 2. 
We take a closed cylinder neighborhood of the disc.
Inside this neighborhood we make the inverse of the 
modification given in Lemma 2. 
This move lowers the level of the 2-handle. This lowering can be continued
until the 2-handle is moved to the level $c$. Thus, moving 2-handles down
is possible if there is an embedded collar that does not meet any 
critical points. 

We have on the level $c$ the 2-handle and on the level $f(p_1)$ the 1-handle.
The boundary mapping is trivial. These two handles cancel each other, i.e.,
the 2-handle fills the hole created by the 1-handle. 
We can construct a new Morse function where these handles do not appear.
The new Morse function $f'$ satisfies Lemma 10
and the genus of $\partial {M^a}$ for $f'$ is $g-1$. $\eop$   

The first generally accepted proof of the Poincar\'e conjecture was given 
by Grigory Perelman in the year 2002. The proof applies advanced techniques.
We have obtained another proof by using elementary methods only:

\th Theorem 1. Every simply connected closed 3-manifold
is homeomorphic to the 3-sphere.

\proof
We may assume, that $M$ is differentiable since if the claim
holds for differentiable 3-manifolds it holds for all, see e.g. [1]. 
Assume, that $M$ is an exotic 
homotopy sphere. Let $f:M\to \R$ be a Morse function as in Lemma 10. 
Assume, that there is a minimum genus $g$ for the Heegaard split induced
by a Morse function satisfying Lemma 10. There is a lower bound to $g$  
as there is no nontrivial Heegaard split of genus 1.    
By Lemma 14 $M$ admits a Heegaard split as in Lemma 10 which has genus $g-1$.
This contradiction shows, that $M$ is $S^3$. $\eop$


\begin{thebibliography}{2}
%

\bibitem {hemp}
J. Hempel, 3-manifolds, Princeton University Press, Princeton, N,J.,1976
\bibitem {smale}
S. Smale, Generalized Poincar\'e's Conjecture in dimensions greater
than 4, Annals of Math.,74 (1961), 391-406
\end{thebibliography}
\end{document}